\documentclass[12pt]{article}
\usepackage[margin=1in]{geometry}

\usepackage{cite}
\usepackage{amsmath,amssymb,amsfonts,mathrsfs}
\usepackage{algorithmic}
\usepackage{graphicx}
\usepackage{textcomp}
\DeclareMathOperator*{\argmax}{argmax}
\DeclareMathOperator*{\argmin}{argmin}

\usepackage[justification=centering]{caption}
\newtheorem{dfn}{Definition}
\newtheorem{thm}{Theorem}
\newtheorem{exm}{Example}
\newtheorem{prp}{Proposition}
\newtheorem{rem}{Remark}
\newtheorem{lmm}{Lemma}

\newtheorem{cor}{Corollary}

\begin{document}

\title{\LARGE \bf
Stabilizing Queuing Networks with Model Data-Independent Control
}
\author{Qian Xie and Li Jin
\thanks{Q. Xie is with the Tandon School of Engineering, New York University, USA. L. Jin is with the UM Joint Institute and the School of Electronic Information and Electrical Engineering, Shanghai Jiao Tong University, China and with the Tandon School of Engineering, New York University, USA (emails: qianxie@nyu.edu, li.jin@sjtu.edu.cn).}}
\maketitle

\thispagestyle{plain}
\pagestyle{plain}

\begin{abstract}
Classical queuing network control strategies typically rely on accurate knowledge of model data, i.e. arrival and service rates.
However, such data are not always available and may be time-variant.
To address this challenge, we consider a class of model data-independent (MDI) control policies that only rely on traffic state observation and network topology.
Specifically, we focus on the MDI control policies that can stabilize multi-class Markovian queuing networks under centralized and decentralized policies.
Control actions include routing, sequencing, and holding.
By expanding the routes and constructing piecewise-linear test functions, we derive an easy-to-use criterion to check the stability of a multi-class network under a given MDI control policy.
For stabilizable multi-class networks, we show that a centralized, stabilizing MDI control policy exists.
For stabilizable single-class networks, we further show that a decentralized, stabilizing MDI control policy exists.
In addition, for both settings, we construct explicit policies that attain maximal throughput and present numerical examples to illustrate the results.
\end{abstract}

\section{Introduction}


Control on multi-class queuing networks has been studied in numerous contexts of transportation, logistics, and communication systems \cite{kumar1995stability,meyn2001sequencing,smith2010dynamic,zhang2018analysis}.
Most existing analysis and design approaches rely on full knowledge of model data, i.e., arrival and service rates, to ensure stability and/or optimality \cite{bertsimas1994optimization}.
However, in many practical settings, such data may be unavailable or hard to estimate, and may be varying over time.
Such challenges motivate the idea of \emph{model data-independent (MDI)} control policies.
MDI control policies select control actions, including routing, sequencing, and/or holding, according to state observation and network topology but independent of arrival/service rates.
Such policies are easy to implement and, if appropriately designed, can resist modeling error or non-stationary environment.
However, the stability of general open multi-class queuing networks with centralized or decentralized MDI control policies has not been well studied.

In this paper, we consider the stability of multi-class queuing networks with throughput-maximizing MDI control policies.
Particularly, we focus on acyclic open queuing networks with Poisson arrivals and exponential service times.
Jobs (customers) are classified according to their origin-destination (OD) information.
Service rates are independent of job classes.
A network is stabilizable if there exists a control policy that ensures positive Harris recurrence of the queuing process, whether the network is open-loop or closed-loop, centralized or decentralized \cite{down1997piecewise}.
By standard results on Jackson networks, stabilizability is equivalent to the existence of a (typically model data-dependent) stabilizing Bernoulli routing policy \cite{gallager2013stochastic}.
We assume that the class-specific arrival rates and the server-specific service rates are unknown to the controller.
The main results are as follows:
\begin{enumerate}
    \item An easy-to-use criterion to check the stability of a multi-class network under a given MDI control policy (Proposition~\ref{prp_stability}).
    \item For a multi-class network, a stabilizing centralized MDI control policy exists if and only if the network is stabilizable (Theorem~\ref{thm_centralized}).
    \item For a single-class  network, a stabilizing decentralized MDI control policy exists if and only if the network is stabilizable (Theorem~\ref{thm_decentralized2}).
\end{enumerate}

Previous works on stability of queuing networks are typically based on full knowledge of model data \cite{dai1995positive,kumar1995stability,foss1998stability,tang2020analysis,sarikaya2011dynamic,dube2009bertrand,savla2011dynamical}.
So far, the best-studied MDI control policy is the join-the-shortest-queue (JSQ) routing policy for parallel queues \cite{foschini1978basic,ephremides1980simple,vvedenskaya1996queueing,eschenfeldt2018join,gupta2007analysis,mukhopadhyay2015analysis,mehdian2017join,tang2020security} or simple networks \cite{bramson2011stability}, which requires only the queue lengths and does not rely on model data \cite{dai2007stability}. 
When and only when the network is stabilizable, i.e., the demand is less than capacity (service rate), the JSQ policy guarantees the stability of parallel queues/simple networks \cite{foley2001join,bramson2011stability} and the optimality of homogeneous servers \cite{ephremides1980simple}.
However, JSQ routing does not guarantee stability of more complex networks \cite{dai2007stability}.
MDI routing for general networks has been numerically evaluated \cite{ling2010global}, but no structural results are available.
Besides, most studies on MDI routing for general networks are not aimed for stability \cite{ren2011traffic,papadimitriou1994complexity,towsley1980queuing,kelly1993dynamic}.
In addition, decentralized dynamic routing has been considered for single origin-destination networks \cite{sarachik1982decentralized,gregoire2014capacity} but not in MDI settings.

To design stabilizing MDI control policies, we first develop a stability criterion (Proposition~\ref{prp_stability}) based on route expansion for queuing networks and explicit construction of a piecewise-linear test function.
The expanded network is essentially a parallel connection of all routes from the set of origins to the set of destinations.
With this expansion, we use insights on the behavior of parallel queues and of tandem queues to construct the test function and derive the stability criterion.
The test function can be used to obtain a smooth Lyapunov function verifying a negative drift condition.
The piecewise-linear test function technique was proposed by Down and Meyn \cite{down1997piecewise}; however, their implementation relies on linear programming formulations to determine parameters of the test function, which depends on model data.
We will extend this technique to the MDI setting using explicitly constructed test functions.

Based on the stability criterion, we design control policies in centralized and decentralized settings.
First, for multi-class networks, we present a stabilizing centralized MDI control policy requiring dynamic routing and preemptive sequencing named JSR policy (Theorem~\ref{thm_centralized}).
The control policy is obtained by minimizing the mean drift of the piecewise-linear test function, and the mean drift is guaranteed to be negative if and only if the network is stabilizable.
The JSR policy, which is centralized and MDI, maximizes throughput among all control policies. Compared with other centralized policies, it does not require knowledge of model data, and compared with other MDI policies (e.g., JSQ), it guarantees stability for any stabilizable networks.
Second, for single-class networks, we present a decentralized routing and holding policy that guarantees stability (Theorem~\ref{thm_decentralized2}).
Such policies can also maximize the throughput since the stabilizability of the network implies that the throughput can be as large as close to the capacity. 
The results are closely related to the theory on the classical JSQ routing policy \cite{dai2007stability} and the decentralized max-pressure control policy \cite{varaiya2013max}.

The rest of this paper is organized as follows.
Section~\ref{sec_model} defines the multi-class queuing network model.
Section~\ref{sec_stability} presents the stability criterion based on route expansion and piecewise-linear test function.
Section~\ref{sec_centralized} and Section~\ref{sec_decentralized} consider the control design problem in centralized and decentralized settings respectively.
Section~\ref{sec_conclude} gives concluding remarks. 
\section{Multi-class queuing network}
\label{sec_model}
Consider an acyclic network of queuing servers with infinite buffer spaces.
Let $\mathcal N$ be the set of \emph{servers}. Each server $n$ has an exponential \emph{service rate} $\bar\mu_n$.
The network has a set $\mathcal S$ of \emph{origins} and a set $\mathcal T$ of \emph{destinations}.
Jobs are classified according to their origins and destinations.
That is, we can use an origin-destination (OD) pair $(S,T)\in\mathcal C$ to denote a \emph{job class}, or simply \emph{class}.
For notational convenience, classes (OD pairs) are indexed by $c=(S_c,T_c)$.
Jobs of class $c$ arrive at $S_c$ according to a Poisson process of rate $\lambda_c\ge0$.
We assume that service rates are independent of job class.

The \emph{topology} of the network is characterized by \emph{routes} between origins and destinations.
We use $|r|$ to denote the number of servers on route $r$.
Let $\mathcal R_c$ be the set of routes between $S_c$ and $T_c$, and define $\mathcal R=\bigcup_{c\in\mathcal C}\mathcal R_c$.
Below is an example network to illustrate the notations.
\begin{exm}
Consider the Wheatstone bridge network in Fig.~\ref{fig_bridge2}.
\begin{figure}[hbt!]
  \centering
  \includegraphics[width=0.4\textwidth]{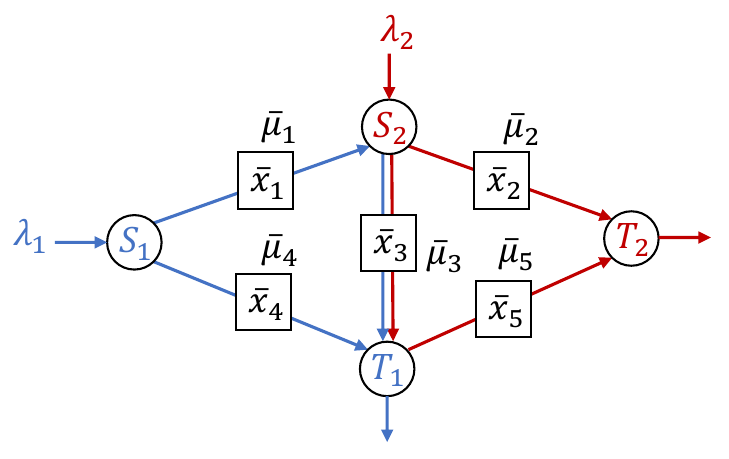}
  \caption{A two-class queuing network.}
  \label{fig_bridge2}
\end{figure}
Two classes of jobs arrive at $S_1$ (resp. $S_2$) with $\lambda_1>0$ (resp. $\lambda_2>0$). The set of servers is $\mathcal N=\{1,2,\ldots,5\}$ and the set of OD-specific routes are
\begin{align*}
    \mathcal R_{1}=\{(1,3),(4)\},\ 
    \mathcal R_{2}=\{(2),(3,5)\}.
\end{align*}
\end{exm}

The \emph{state} of the network is defined as follows.
Let $\bar x=[\bar x_n^c]_{n\in\mathcal N,c\in\mathcal C}$ be the vector of class-specific \emph{job numbers}, where $\bar x_n^c$ is the number of jobs of class $c$ in server $n$, either waiting or being served.
Let $\bar{\mathcal X}$ be the space of $\bar x$.
We use $\bar X(t)$ to denote the state of the queuing process at time $t$.%

We consider three types of \emph{control actions}, viz. {routing}, {sequencing}, and {holding}.
All control actions are essentially Markovian (in terms of $\bar{x}$ plus additional auxiliary states) and are applied at the instant of \emph{transitions}, which include the arrival of a job at an origin or the completion of service at a server.
\emph{Routing} refers to allocating an incoming job to a server downstream to the origin or allocating a job discharged by a server to another downstream server.
\emph{Sequencing} refers to selecting a job from the waiting queue to serve.
The default sequencing policy is the first-come-first-serve (FCFS) policy.
For the multi-class setting, we consider the preemptive-priority that can terminate an ongoing service and start serving jobs from another class, while the job with incomplete service is sent back to the queue.
\emph{Holding} refers to holding a job that has completed its service in the server while blocking the other jobs in the queue from accessing the server.

Following \cite{dai1995stability}, we say that a queuing network is \emph{stable} if the queuing process is positive Harris recurrent. 
For details about the notion of positive Harris recurrence for queuing networks, see \cite{dai1995positive,dai1995stability,down1997piecewise}.
Finally, we say that the network is \emph{stabilizable} if a stabilizing control exists.
One can check the stabilizability using the following result:
\begin{lmm}
\label{lmm_stabilizability}
An open acyclic queuing network is stabilizable if and only if there exists a vector $[\xi_r]_{r\in\mathcal R}$ such that
\begin{align*}
    &\xi_r\ge0,\quad\forall r\in\mathcal R,\\
    &\lambda_c=\sum_{r\in\mathcal R_c}\xi_r,\quad\forall c\in\mathcal C,\\
    &\sum_{r\in\mathcal R:n\in r}\xi_r<\bar\mu_n,\quad\forall n\in\mathcal N.
\end{align*}
\end{lmm}
The proof and implementation are straightforward. 
\section{Stability criterion}
\label{sec_stability}

In this section, we derive a stability criterion for multi-class networks under given control policies.
The techniques that we use include the route expansion of the original network and the explicit construction of a piecewise-linear test function based on the network topology.
In Section~\ref{sub_route}, we construct an expanded network based on the original network.
In Section~\ref{sub_stability}, we apply a piecewise-linear test function to the expanded network to obtain a stability criterion (Proposition~\ref{prp_stability}) for both the expanded and the original networks.

\subsection{Route expansion}
\label{sub_route}

\begin{figure}[hbt!]
  \centering
  \includegraphics[width=0.45\textwidth]{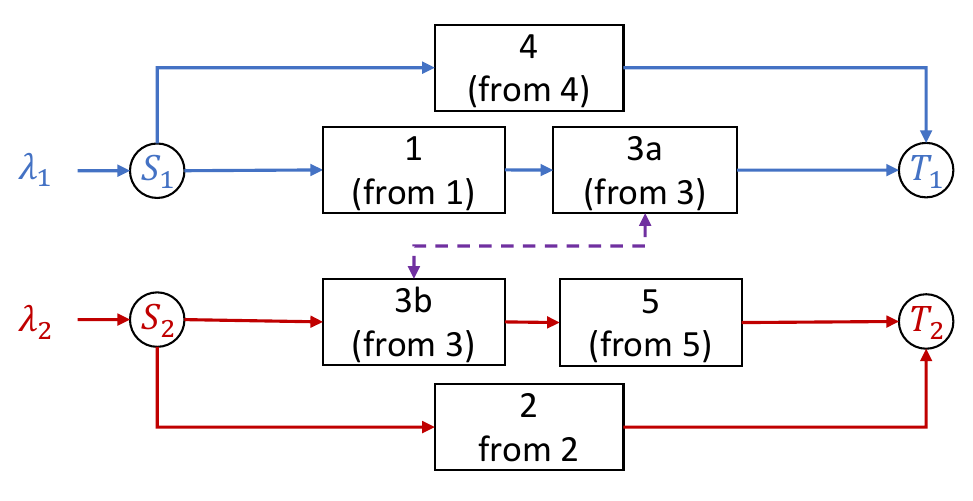}
  \caption{Route expansion of the network in Fig.~\ref{fig_bridge2}. }
  \label{fig_exp2}
\end{figure}

For the convenience of constructing test function, we first introduce the route expansion. \emph{Route expansion} refers to the construction of an \emph{expanded network} based on the topology of \emph{original network} (defined in Section~\ref{sec_model}).
The high-level idea is to decompose the network into routes, and the specific procedures are:
\begin{enumerate}
    \item Place all routes $\mathcal R$ in the original network in parallel.
    \item Add two-way connections between duplicates of servers in the original network.
\end{enumerate}
For example, Fig.~\ref{fig_exp2} shows the expanded network constructed from the original network in Fig.~\ref{fig_bridge2}. 

We call ``servers" in the expanded network as \emph{subservers}, since they are obtained by duplicating actual servers in the original network.
Subservers are indexed by $k$, $c_k\in\mathcal C$ is the class index, $r_k\in\mathcal R$ is the route index, and $i_k\in\{1,2,\ldots,|r_k|\}$ is the numbering of subserver $k$ on route $r_k$.
We use $k \in r$ to refer to that subserver $k$ is on route $r$.
Let $\mathcal K$ be the set of all subservers and $\mathcal K_c$ be the set of subservers with $c_k=c$.
We use $n_k\in\mathcal N$ to denote the actual server that corresponds to subserver $k$. In addition, let $k_p$ (resp. $k_s$) denote the subserver immediately upstream (resp. downstream) to subserver $k$.

The \emph{state} of the expanded network is $x=\{x_{k};k\in\mathcal K\}$, denoting the vector of number of jobs in subserver $k$. 
The expanded state space is $\mathcal{X}=\mathbb Z_{\ge0}^{|\mathcal K|}$. 
Note that the states of the expanded network and the states of the original network are related by
\begin{align}
    &\bar x^c_n=\sum_{k\in\mathcal K:n_k=n, c_k=c}x_{k},
    \quad k\in\mathcal K,\label{eq_xbarx}
\end{align}
for each $n\in\mathcal N$.

The routing policy is characterized by $\pi:\mathcal{X}\to[0,1]^{|\mathcal{K}|^2}$, where $\pi_{k,k'}$ is the probability that a job is routed from subserver $k$ to subserver $k'$.

The holding policy is characterized by $\zeta:\mathcal{X}\to\{0,1\}^{|\mathcal{K}|}$, where $\zeta_k$ specifies whether subserver $k$ is holding ($\zeta_k(x)=0$) or not holding ($\zeta_k(x)=1$) when the current state is $x$.

Two subservers $k$ and $k'$ are \emph{duplicating} if $n_k=n_{k'}$.
Note that the service rates of duplicating subservers are coupled in the sense that for each server $n\in\mathcal N$, at a given time, at most one subserver $k$ such that $n_k=n$ can be actively serving jobs, or \emph{active}.
This can be modeled as an \emph{imaginary service rate control policy} $\mu:\mathcal{X}\to\mathbb{R}^{|\mathcal{K}|}$ such that the service rate $\mu_{k}(x)$ of subserver $k$ satisfies
\begin{align*}
    \sum_{k:n_k=n}\mu_{k}(x)\le\bar\mu_n,
    \quad\forall x\in\mathcal{X}.
\end{align*}
Such control policy is essentially equivalent to the class-based preemptive sequencing in the original network.

Note that $\{X(t): t\geq0\}$ is a Markov process, and the positive Harris recurrence refers to that there exists a unique invariant measure $\nu$ on $\mathcal{X}$ such that for every measurable set $D\subseteq\mathcal{X}$ with $\nu(D)>0$ and for every initial condition $x\in\mathcal{X}$,
$$
\Pr\{\tau_D<\infty|X(0)=x\}=1,
$$
where $\tau_D=\inf\{t\ge0:X(t)\in D\}$. Also, though $\{\bar{X}(t): t\geq0\}$ is not a Markov process, it will eventually converge to a steady state distribution.

The route expansion technique not only expands the network but also decomposes the state variables.
Jobs can move along the expanded network using two transition mechanisms.
One is \emph{actual transition}, referring to moving a job from subserver $k$ (or an origin) to its downstream subserver $k_s$ (or a destination).
{The other is \emph{imaginary transition} that moves a job from one subserver $k$ to a duplicating subserver $k'$ thereof, see imaginary switch in Section~\ref{sec_decentralized}.}
Imaginary transitions always occur instantaneously.
Note that an actual transition corresponds to a transition in the original network, while an imaginary transition does not; this is also revealed in \eqref{eq_xbarx}.

One can always map a control action in the expanded network to the original network.
However, an MDI control policy may not exist on the state space of the original network; we do need an expanded state space for MDI control.
In addition, we allow {imaginary} control actions in the expanded network, including \emph{imaginary service rate control} and \emph{imaginary switch}; see Section~\ref{sec_centralized} and Section~\ref{sec_decentralized}.
Such imaginary actions only make sense in the expanded network and do not correspond to actual service rate control or switch in the original network.

\subsection{Stability of the expanded network}
\label{sub_stability}
After introducing the expanded network and the mathematical definition of the control policy, the main question of this paper can be expressed in a formal way as follows.

Given an expanded queuing network and a control policy $\phi=(\pi,\mu,\zeta)$, how can we tell if the control policy $\phi$ is stabilizing, i.e., the expanded network is stable under $\phi$?

The answer of this question will be given in Proposition~\ref{prp_stability}. Before that, we need to introduce the test function technique first. As opposed to linear programming-based construction in \cite{down1997piecewise}, we provide an explicit construction, where parameters of the test function do not rely on solving any optimization problems. The high level idea is to identify the bottlenecks and their upstream subservers. Our construction is based on the route expansion described in the previous subsection.
\begin{enumerate}
\item For each class $c\in\mathcal C$ and expanded state $x\in\mathcal X$, define
\begin{align*}
    g_c(x):=\max_{\substack{K_c\subseteq\mathcal K_c:\\\kappa\in K_c\Rightarrow \kappa_p\in K_c}}\sum_{k\in K_c}a_kx_k,
\end{align*}
where $a_k\in(0,1)$ is a parameter.

\item 
Define a piecewise-linear \emph{test function}
\begin{align*}
    V(x):=\max_{C\subseteq\mathcal C}\sum\limits_{c\in C}b_cg_c(x),
\end{align*}
where $b_c\in(0,1)$ is a parameter.
\end{enumerate}


We call $V(x)$ the test function rather than the Lyapunov function, since strictly speaking, a smooth Lyapunov function should be developed based on the piecewise-linear test function to verify the Foster-Lyapunov stability criterion.
Down and Meyn \cite{down1997piecewise} showed that as long as a piecewise-linear test function can be determined, one can always smooth it to obtain a qualified $C^2$ Lyapunov function.

\begin{rem}
 The test functions we proposed in this work are MDI. But generally speaking, they do not need to be MDI since it does not affect the control policies to be MDI.
\end{rem}

\begin{dfn}[Dominance]
\label{dfn_dominance}
Consider state $x\in\mathcal X$. 

\begin{enumerate}
    \item We call $C^*$ a set of \emph{dominant} classes if
    $$
     C^*\in\argmax_{C\subseteq\mathcal C}\sum_{c\in C}b_{c}g_{c}(x).
    $$
    Each class $c\in C^*$ is a \emph{dominant} class. 
    \item We call $K^*_c$ a set of \emph{dominant} class-c subservers if
    $$K^*_c\in\argmax\limits_{\substack{K_c\subseteq\mathcal K_c:\\k\in K_c\Rightarrow k_p\in K_c}}\sum_{k\in K_c}a_kx_k.$$
    Each subserver $k\in K^*_c$ is a \emph{dominant} class-c subserver.
    \item A route $r\in\mathcal R_c$ is \emph{dominant} if it includes dominant class-c subservers, i.e. there exists dominant class-c subserver $k\in K^*_c$ such that $k\in r$.
    
    Let $R_c$ be the set of dominant class-c routes.
    \item A subserver $b\in K^*_c$ is called a \emph{bottleneck} if it is a dominant class-c subserver while its immediate downstream subserver $b_s\notin K^*_c$ is not.
\end{enumerate}
\end{dfn}

\begin{rem}
A route or server is dominant if changes in its traffic state immediately affect the test function $V$.
\end{rem}

    
    
    

A \emph{regime} $X$ of the piecewise-linear test function is a subset of $\mathcal X$ such that there exist $C^X\subseteq\mathcal C$, $K^X=\bigcup_{c\in C^X}K_{c}^X\subseteq\mathcal K$, and $R^X=\bigcup_{c\in C^X}R_c^X\subseteq\mathcal R$ where $C^X,K_c^X,R_c^X$ are dominant for each $x\in X$, i.e., the test function is linear over $X$.
Let $\mathscr X$ be the set of regimes; note that $\bigcup_{X\in\mathscr X}X=\mathcal X$.

\begin{dfn}[Mean velocity and drift]
Consider a multi-class network with state $x\in\mathcal{X}$ under an expanded control policy $\phi=(\pi,\mu,\zeta)$.
\begin{enumerate}
    \item The \emph{mean velocity} at state $x$ is a function $v:\mathcal{X}\to\mathbb R^{|\mathcal K|}$ such that for each $k\in\mathcal K$,
\begin{align*}
    v_{k}(x):=\sum_{c\in\mathcal C}\lambda_{c}\pi_{S_{c},k}^c(x)
    +\mu_{k_p}(x)\zeta_{k_p}(x)
    -\mu_{k}(x)\zeta_k(x).
\end{align*}
where $\pi^c_{S_c,k}$ is the probability that a class-$c$ job is routed from origin $S_c$ to subserver $k$, while $\mu_k$ and $\zeta_k$ are the controlled service rate and the holding status of the subserver $k$ respectively.

\item Given $X\in\mathscr X$ such that $x\in X$, the \emph{mean drift over $X$} is given by
\begin{align*}
    &D^X(x):=\sum_{c\in C^X}b_c\sum_{k\in K^X_{c}}a_kv_{k}(x).
\end{align*}
\end{enumerate}
\end{dfn}

\begin{rem}
In our subsequent analysis, the mean drift $D^X(x)$ of the test function will play the role of infinitesimal generator applied to a Lyapunov function; see \cite{down1997piecewise} for the connection between the test function and the Lyapunov function.
\end{rem}

The main result of this section is as follows:
\begin{prp}
\label{prp_stability}
Consider a multi-class network under the expanded control policy $\phi$. Suppose there exist constants $M<\infty$, $\epsilon>0$, and $a_k, b_c\in(0,1)$ ($\forall$ $c\in\mathcal C$, $k\in K_c$), such that for each $X\in\mathscr X$ and each $x\in X$ where $|x|=\sum\limits_{k\in\mathcal K}x_k> M$,
\begin{align}
    \sum_{c\in C^X}b_c\sum_{k\in K^X_{c}}a_kv_{k}(x)\le-\epsilon.
    \label{eq_betaalpha}
\end{align}
Then, the network is stable.
\end{prp}

\noindent\emph{Proof.}
Consider the following test function $V(x)$:
$$
V(x)=\sum_{c\in C^X}b_c\sum_{k\in K^X_{c}}a_kx_{k}.
$$
By \eqref{eq_betaalpha}, the mean drift satisfies
\begin{align*}
    D^X(x)&=\sum_{c\in C^X}b_c\sum_{k\in K^X_{c}}a_kv_{k}(x)\le-\epsilon
    \quad x:|x|>M.
\end{align*}
One can then apply \cite[Theorem 1]{down1997piecewise} and \cite[Lemma 5]{down1997piecewise} to obtain the stability of the network.
\hfill$\square$

As a benchmark, the approach in \cite[Theorem 1]{down1997piecewise} requires solving linear programs to obtain parameters of the test functions in Proposition~\ref{prp_stability}, while our approach explicitly constructs the parameters (see Section \ref{sec_centralized} and Section \ref{sec_decentralized}). Moreover, the proposed control, which is independent of model data, guarantees stability if and only if the network is stabilizable (see Theorem~\ref{thm_centralized} and Theorem~\ref{thm_decentralized2}), while the approach in \cite{down1997piecewise} relies on knowledge of model data. 
\section{Centralized control for multiple classes}
\label{sec_centralized}

In this section, we consider the ``join-the-shortest-route (JSR)'' policy (a joint routing and sequencing policy) for centralized control. The JSR policy is MDI and constructed based on the expanded network. We will show that it is stabilizing if and only if the network is stabilizable. 

The test functions are constructed as follows.
\begin{enumerate}
    \item For each class $c\in\mathcal C$, each route $r\in\mathcal R_c$, and each expanded state $x\in\mathcal X$, let
    \begin{align*}
        &f_{r}(x):=\max_{k\in r}\alpha^{i_k-1}\sum_{j:i_j\le i_k}x_{j},\\
        &g_c(x):=\max_{R_c\subseteq\mathcal R_c}\beta^{|R_c|-1}\sum_{r\in{R_c}}f_{r}(x),
    \end{align*}
    where $\alpha\in(0,1),\beta\in(0,1)$ are constant parameters.
    \item 
    The piecewise-linear \emph{test function} is given by
    \begin{align*}
        V(x):=\max_{C\subseteq\mathcal C}\gamma^{|C|-1}\sum_{c\in C}g_{c}(x),
    \end{align*}
    where $\gamma\in(0,1)$ is a constant parameter.
\end{enumerate}
Let the parameters be such that
\begin{align}
    \alpha=\beta\geq\frac{|\mathcal R|-1}{|\mathcal R|},\quad
    \gamma\geq\frac{|\mathcal C|-1}{|\mathcal C|},
    \label{eq_parameter}
\end{align}
and follow the notions of dominance accordingly (see Definition~\ref{dfn_dominance}).
Note that such MDI parameters $\alpha,\beta,\gamma$ always exist.
The control that we consider in this subsection only depends on $\alpha,\beta,\gamma$ and is thus MDI.
Specifically, we define the JSR policy as follows:

\begin{dfn}[Join-the-shortest-route (JSR) policy]
\label{dfn_JSR}
$ $

\begin{enumerate}
    \item (Routing) At an origin $S$, an incoming job of class $c$ is allocated to the route $r^*\in\mathcal R_c$ such that
    \begin{align*}
        r^*\in\argmin_{r\in\mathcal R_c}f_{r}(x).
    \end{align*}
    If there is only one minima, then $r^*$ must be a non-dominant route. Otherwise, let $b^*$ be the bottleneck on route $r^*$.
    Then, an incoming job of class $c$ is allocated to the route $r^*\in\mathcal R_c$ with the largest $i_{b^*}$, which is denoted by $i_c$.
    Further ties are randomly broken.
    
    \item (Imaginary service rate control) Let $\mathcal K_n$ be the set of subservers corresponding to server $n$ and let $\mathcal B$ be the set of bottlenecks for a given $x$. Then, a subserver $k\in\mathcal K_n$ is activated if $k\in\mathcal B$. If multiple subservers are in $\mathcal K_n\cap\mathcal B$, then activate the subserver $k^*$ such that
    \begin{align*}
        k^*=\argmin_{k\in\mathcal K_n\cap\mathcal B}\{i_{c_k}+|\mathcal{R}_{c_{k}}|\};
    \end{align*}
    ties are randomly broken. This is to ensure that the bottlenecks are active to discharge jobs and only one of the duplicating subservers can be active.
    
\end{enumerate}
\end{dfn}

The main result of this section is the following:

\begin{thm}[Stability of JSR policy]
\label{thm_centralized}
The JSR policy stabilizes a multi-class network if and only if the network is stabilizable.
\end{thm}

This theorem implies that the JSR policy is also throughput-maximizing, as long as the network is stabilizable, i.e., the demand is less than the total capacity. Note that the stabilizability can be easily checked using Lemma~\ref{lmm_stabilizability}.

In the rest of this section, we apply Theorem~\ref{thm_centralized} to study the stability of the Wheatstone bridge network under the JSR policy (Subsection~\ref{sub_centralized_example}) and then prove this theorem  (Subsection~\ref{sub_centralized_proof}).

\subsection{Numerical Example}
\label{sub_centralized_example}
Consider the network in Fig.~\ref{fig_bridge2} and suppose that $\lambda_1=\lambda_2=\lambda=1$ and $\bar\mu_n=\mu=1$ for $n=1,2,4,5$ and $\bar\mu_3=\frac14$.
This example is for illustrating the route expansion and the test function construction.

Note that under the above model parameters, the decentralized JSQ policy is destabilizing. To see this, $\bar\mu_1=\bar\mu_4$ implies that on average, class-1 jobs are evenly distributed between server 1 and server 4.
Thus, the average departure rate of class-1 jobs from server 1 is $\frac12$, which exceeds the service rate of server 3. Therefore, the queue at server 3 is unstable.
The main reason that the JSQ policy is destabilizing is the ignorance of downstream congestion. As $\bar X_3(t)$ gets large, a reasonable action is to allocate fewer class-1 jobs to server 1. However, the JSQ policy disallows such far-sighted decisions.

An alternative centralized stabilizing routing policy can be the following JSR policy:
\begin{enumerate}
    \item A class-1 job arriving at $S_1$ is routed to server 1 if $\bar X_1^1(t)+\bar X_3^1(t)<\bar X_4^1(t)$, to server 4 if $\bar X_1^1(t)+\bar X_3^1(t)>\bar X_4^1(t)$, and randomly otherwise.
    \item A class-2 job arriving at $S_2$ is routed to server 3 if $\bar X_3^2(t)+\bar X_5^2(t)<\bar X_2^2(t)$, to server 2 if $\bar X_3^2(t)+\bar X_5^2(t)>\bar X_2^2(t)$, and randomly otherwise.
    \item The dominant class has a higher priority.
\end{enumerate}
That is, when jobs are routed at $S_1$, the decision is based on not only the local state ($\bar X_1(t)$ and $\bar X_4(t)$), but also the state further downstream ($\bar X_3(t)$).

The expanded network is shown in Fig.~\ref{fig_exp2}. Each block in the figure represents a subserver.
In particular, subservers 3a and 3b are decomposed from server 3; the other servers are remained.
Solid arrows correspond to actual transitions in an original network, while dashed arrows correspond to imaginary transitions between duplicating subservers.

In the expanded network, a job can move along both solid and dashed arrows.
The color of an arrow shows which class can move along it: blue means class $(S_1,T_1)$, red means $(S_2,T_2)$, and purple means both.
For ease of presentation, we label $(S_1,T_1)$ as class 1 and $(S_2,T_2)$ as class 2.
For example, a job of class $(S_1,T_1)$ can visit subservers 4, 1, 3a, 3b and the destination $T_1$.

In the expanded network, the JSR policy works as follows. 
\begin{enumerate}
    \item A class-1 job arriving at $S_1$ is routed to subserver 4 if $X_{4}(t)<X_{1}(t)+X_{3a}(t)$, to subserver 1 if $X_{4}(t)>X_{1}(t)+X_{3a}(t)$, and randomly otherwise.
    
    \item A class-2 job arriving at $S_2$ is routed to subserver 3b if $X_{3b}(t)+X_{5}(t)<X_{2}(t)$, to subserver 2 if $X_{3b}(t)+X_{5}(t)>X_{2}(t)$, and randomly otherwise.
    
    \item If subserver 3a is dominant while subserver 3b is non-dominant, and server 3 is serving a class-2 job, then server 3 preempts the class-2 job being served in 3b to the class-1 job in 3a, and vice versa. If both subserver 3a and subserver 3b are dominant, then server 3 gives priority to the class-2 job since the index of 3b is smaller.
\end{enumerate}
By Theorem~\ref{thm_centralized}, the network can be stabilized by the JSR policy if and only if
\begin{align*}
    \lambda_1<2,\
    \lambda_2<2,\
    \lambda_1+\lambda_2<\frac94.
\end{align*}
We use the following parameters for the test function:
$$\alpha=\beta=\gamma=\frac34,\  \epsilon=\Big(\frac34\Big)^5.$$
One can verify that the above parameters satisfy \eqref{eq_parameter} and Proposition~\ref{prp_stability} by considering the following cases:
\begin{enumerate}
    \item Only one route is dominant. In this case, an incoming job is always allocated to a non-dominant route, leading to non-positive contribution to the mean drift:
    \begin{align*}
        D^X(x)\le-\gamma\beta\alpha\mu=-\Big(\frac34\Big)^3\le-\epsilon.
    \end{align*}
    \item Two routes with different OD pairs are dominant. This case is analogous to the previous case:
    \begin{align*}
        D^X(x)\le-\gamma\beta\alpha\mu=-\Big(\frac34\Big)^3\le-\epsilon.
    \end{align*}
    \item Two routes with the same OD pair or more than two routes are dominant. In such cases, the mean drift satisfies
    \begin{align*}
        D^X(x)\le\gamma\beta^3(\lambda-\mu-\alpha\mu)=-\Big(\frac34\Big)^5\le-\epsilon.
    \end{align*}
\end{enumerate}
Consequently, the network is stable under the MDI JSR policy.

\subsection{Proof of Theorem~\ref{thm_centralized}}
\label{sub_centralized_proof}
In this subsection, we will show the sufficiency and the necessity respectively, based on the connection between the sign of the mean drift and the stabilizability of the network.
When analyzing the mean drift, we consider two parts: external arrivals and internal transmission.
We first show that any internal transmission does not positively contribute to the mean drift and then show that any positive contribution from external arrivals can always be compensated by internal transmissions.

\subsubsection{Internal transmissions}
\label{sub_internal}
Note that under the JSR policy, every job remains on the route assigned to the job when it enters the network.
Hence, internal transmissions only occur between subservers on the same route.

Given $x$, consider an internal transmission from subserver $k$ to subserver $j$; this implicitly requires $x_k\ge1$.
The definition of dominance ensures that if $j$ is dominant, then so is $k$.
Hence, we need to consider the following cases:
\begin{enumerate}
    \item If $k$ and $j$ are both dominant, the transmission leads to zero contribution to the mean drift $D^X(x)$ for all $X$ such that $x\in X$.
    \item If $k$ is dominant and $j$ is non-dominant, the transmission leads to the following contribution to the mean drift:
    \begin{align*}
    -\alpha^{i_k-1}\mu_{k}(x)\le0.
    \end{align*}
\end{enumerate}
Hence, internal transmissions never lead to positive contribution to the mean drift.

\subsubsection{External arrivals}
\label{sub_external}

Given $x\ne0$, consider a regime $X\in\mathscr X$ such that $x\in X$.
For each $c\in\mathcal{C}$, the JSR policy ensures that if there exists a non-dominant route in $\mathcal R_c$, then an incoming job must be allocated to a non-dominant route in $\mathcal R_c$, leading to non-positive contribution to the mean drift. 
Hence, we only need to consider dominant classes $c$ such that every route in $\mathcal R_c$ is dominant, i.e. $R^X_c=\mathcal R_c$.
Recall that $C^{*}\subseteq\mathcal{C}$ is the set of dominant classes.
The part of the mean drift associated with $c\in C^{*}$ satisfies
\begin{align*}
    D^X_c(x)&\le\gamma^{|C^*|-1}\beta^{|\mathcal R_c|-1}\Big(\alpha^{i_c-1}\lambda_c-\sum_{b\in\mathcal B^X:c_b=c}\alpha^{i_b-1}\mu_{b}(x)\Big)\\
    &:=\gamma^{|C^*|-1}\Delta_c^X(x)
\end{align*}
over any regimes of the piecewise-linear test function, where $i_c$ is given in Definition~\ref{dfn_JSR}.

\begin{lmm}
\label{lmm_bottleneck}
When $x\ne0$, there is no empty bottleneck, i.e.
\begin{align}
    x_{b}\ge1.
    \label{eq_xri>1}
\end{align}
\end{lmm}

\noindent\emph{Proof.}

Since $x\ne0$ and $r_b$ is dominant, we have
$$
    \sum_{k:i_k\le i_b}x_{k}>0.
$$
If $i_b=1$, then the above inequality directly implies \eqref{eq_xri>1}.

Now consider the case that $i_b\ge 2$.
Since $b$ is a bottleneck, we have
\begin{align*}
    \alpha^{i_b-1}\sum_{k\in r_b:i_k\le i_b}x_{k}
    \ge\alpha^{i_b-2}\sum_{k\in r_b:i_k\le i_b-1}x_{k},
\end{align*}
which implies
\begin{align*}
    x_{b}\ge(1-\alpha)\sum_{k:i_k\le i_b}x_{k}>0
\end{align*}
and thus we have \eqref{eq_xri>1}.
\hfill$\square$

Lemma \ref{lmm_bottleneck} is to ensure that the bottlenecks are none-empty to discharge jobs and thus contribute negative terms to the drift.




Next, we show the sufficiency of Theorem~\ref{thm_centralized}. Based on the definition of the routing policy (see Definition \ref{dfn_JSR}), $\forall b\in\mathcal B^X$, we have $i_b\leq i_c$ when the incoming job is allocated to a dominant route. Then
\begin{align*}
    \sum_{c\in{C^*}}\Delta_c^X(x)
    &\le\sum_{c\in{C^*}}\beta^{|\mathcal R_c|-1}\alpha^{i_c-1}\Big(\lambda_c-\sum_{b\in \mathcal B^X:c_b=c}\mu_{b}(x)\Big)\\
    &=\sum_{c\in{C^*}}\alpha^{i_c+|\mathcal R_c|-2}\Big(\lambda_c-\sum_{b\in \mathcal B^X:c_b=c}\mu_{b}(x)\Big)\\
    &:=\sum_{c\in{C^*}}\alpha_c\beta_c
\end{align*}
Without loss of generality assume $C^*=\{1,2,\cdots,m\}$ and
$$i_1+|\mathcal{R}_1|\le i_2+|\mathcal{R}_2|\le\cdots\le i_m+|\mathcal{R}_m|.$$

Then by using Abel transformation (summation by parts), the right hand side of the above inequality (abbr. RHS):
\begin{align*}
    RHS = \sum_{i=1}^{m-1}(\alpha_i-\alpha_{i+1})\sum_{j=1}^{i}\beta_j+\alpha_m\sum_{j=1}^{m}\beta_m.
\end{align*}
Based on the assumption, we have $\alpha_i\geq\alpha_{i+1}$ and
\begin{align*}
    \sum_{j=1}^{i}\beta_{j}&=\sum_{j=1}^{i}\Big(\lambda_{j}-\sum_{b\in \mathcal B^X:c_b=j}\mu_{b}(x)\Big)\\
    &=\sum_{j=1}^{i}\lambda_j-\sum_{n_b:b\in\mathcal B_i^X}\bar\mu_{n_b}\\
    &=\sum_{j=1}^{i}\lambda_j-\sum_{n\in\mathcal{N}_i}\bar\mu_{n}\\
    &<0,
\end{align*}
where $\mathcal{B}_i$ is the set of bottlenecks in the first $i$ classes and $\mathcal{N}_i$ is the min-cut of the original network with the first $i$ classes. Here we use the definition of the imaginary service rate control (see Definition~\ref{dfn_JSR}) and Lemma~\ref{lmm_stabilizability}. 

Since $RHS<0$, we have $\sum_{c\in{C^*}}\Delta_c^X(x)<0$ and thus $D_c^X(x)<0$. Then by noting that internal transmissions lead to non-positive contributions to the mean drift, we have
\begin{align*}
    D^X(x)\le\sum_{c\in{C^*}}D^X_c(x)
    <0,
\end{align*}
which implies stability.
\hfill$\square$

Finally, the necessity is apparent: if a network is not stabilizable, then there exists no MDI control that can stabilize the network. 
\section{Decentralized control for a single class}
\label{sec_decentralized}

For a single-class network, we can drop the class index and use $x_k$ to denote the number of jobs in subserver $k$.
Note that such network has a single origin and a single destination. Again we can do route expansion on such network. 

We consider a decentralized MDI control policy as follows.
\begin{dfn}[JSQ with artificial spillback]
The JSQ with artificial spillback (JSQ-AS) policy is as follows:
\begin{enumerate}
    \item (Routing) A discharged job is routed to the shortest downstream queue, with ties randomly broken.
    
    \item (Holding) For each subserver $k$, any job which has finished the service will be held if and only if $X_{k_s}(t)\ge X_{k}(t)$.
    
    
    \item (Imaginary switch) When a dominant subserver $k$ is inactive while its non-dominant duplicate $k'$ is active, and both are not in the holding status, then the job in $k'$ is moved to (and discharged from) $k$ after service, and then routed to the downstream of $k$ (i.e., $k_s$).
\end{enumerate}
\end{dfn}

Note that under the holding policy, the process $\{X(t);t\ge0\}$ admits an invariant set $\mathcal Q\subseteq\mathcal X$ given by
\begin{align}
    \mathcal Q:=\{x\in\mathcal X:x_{k_s}\le x_{k}, k\in\mathcal K\}.
\label{invariant}
\end{align}

Since we consider the long-time stability of the network, it suffices to consider the states in an invariant set.
The above result indicates that in the invariant set $\mathcal Q$, the queue size of any subserver is upper-bounded by the queue size of its immediate upstream subserver.

The JSQ-AS policy is decentralized in the sense that control actions on subserver $k$ only depend on local traffic information: the number of jobs in duplicate subservers $\{x_{k'}:n_k=n_{k'}\}$ and that in immediate downstream subservers $\{x_{{k'}_s}:n_k=n_{k'}\}$.
A key characteristic of such policies is that congestion information can propagate through the network via the forced holding: if a subserver becomes congested (i.e. $x_{k}$ gets large), the congestion will propagate to the upstream subservers in a cascading manner (``artificial spillback'').
Importantly, such artificial spillback does not undermine throughput like the natural spillback caused by the limited buffer size.
The reason is that though congestion can propagate, the queue size in any downstream subserver is not upper-bounded.
Artificial spillback is the main difference between the JSQ-AS policy and the classic JSQ policies.

Note that though the JSQ-AS policy is constructed based on the expanded network, its actions can always be converted to the ones in the original network.
Importantly, the decentralized control in the expanded network must also be decentralized in the original network.
Also note that the imaginary switch has no impact on the original network or the test function.

The main result of this section is as follows:

\begin{thm}[Stability of JSQ-AS policy]
\label{thm_decentralized2}
For the route expansion of a single-class network, the JSQ-AS policy is stabilizing if and only if
\begin{align}
    \lambda<\bar\mu^{mc},
    \label{eq_lambda<mu}
\end{align}
where $\bar\mu^{mc}$ is the min-cut service rate of the original network.
\end{thm}

This theorem implies that JSQ-AS policy is also a throughput-maximizing policy since we allow any throughput that satisfies \eqref{eq_lambda<mu}.

In the rest of this section, we apply Theorem~\ref{thm_decentralized2} to study the stability of the Wheatstone bridge network under the JSQ-AS policy (Subsection~\ref{sub_decentralized_example}) and then prove this theorem  (Subsection~\ref{sub_decentralized_proof}).

\subsection{Numerical Example}
\label{sub_decentralized_example}
\begin{figure}[hbt!]
  \centering
  \includegraphics[width=0.48\textwidth]{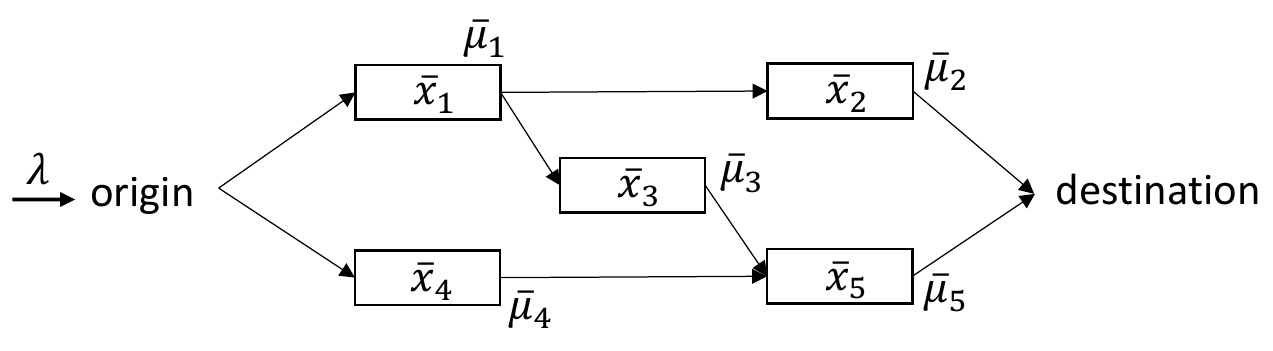}
  \includegraphics[width=0.48\textwidth]{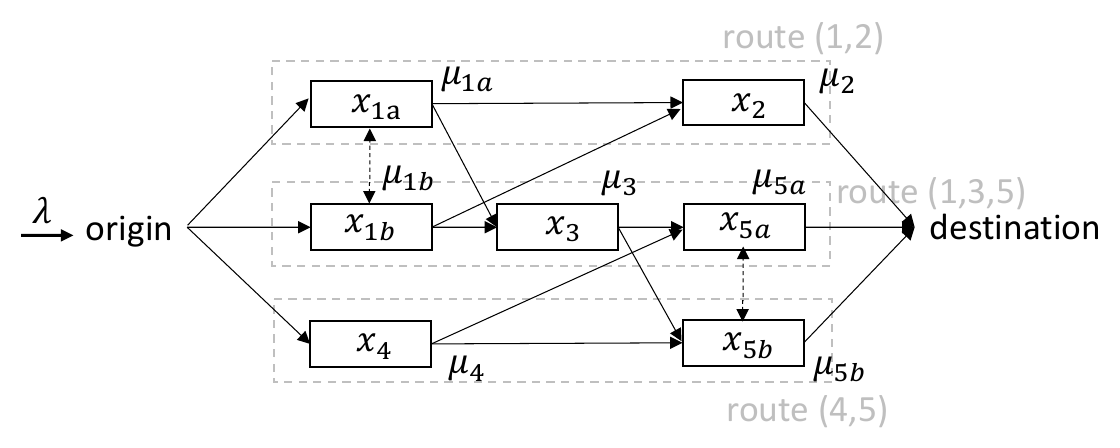}
  \caption{A single-class queuing network and its expanded network.}
  \label{fig_single}
\end{figure}
Consider the original network with route expansion in Fig.~\ref{fig_single}. Again suppose that $\lambda=1$, $\bar\mu_n=\frac34$ for $n=1,2,\cdots,5$. Similarly, with the above parameters, the JSQ policy is destabilizing since the queue at server 5 is unstable. However, in the decentralized setting, the control actions can only depend on the local state, say the routing decision at the origin can be based on $\bar X_1(t)$ and $\bar X_4(t)$, but not $\bar X_5(t)$. A remedy is to introduce the holding policy (artificial spillback) to the JSQ policy so that the downstream congestion can be relieved and the local state can somehow reflect the states further downstream.

In the expanded network, server 1 is decomposed into subserver 1a and 1b, server 5 is decomposed into subserver 5a and 5b. The states in the original network and those in the expanded network satisfy $\bar X_1(t)=X_{1a}(t)+X_{1b}(t)$ and $\bar X_5(t)=X_{5a}(t)+X_{5b}(t)$. The initial states of the expanded network can be not unique. Say the initial queue size of server 1 is 2, then the initial queue sizes of subserver 1a and 1b can be 2, 0 or 1, 1 or 0, 2 respectively. Then the states are updated based on the model and our JSQ-AS policy. For example, the routing decision at the origin is based on $X_{1a}(t)$, $X_{1b}(t)$ and $X_4(t)$ rather than $\bar X_1(t)$ and $\bar X_4(t)$; a job which has just finished the service at server 3 will be held if $X_{5a}(t)\geq X_3(t)$, once released, it will be routed to the shorter downstream queue by comparing $X_{5a}(t)$ and $X_{5b}(t)$.



\subsection{Proof of Theorem~\ref{thm_decentralized2}}
\label{sub_decentralized_proof}
This proof uses the connection between the stabilizability condition \eqref{eq_lambda<mu} and the sign of the mean drift. But before showing the mean drift is negative, we first present the explicit MDI piecewise-linear test function and several key lemmas that can help us analyze the mean drift.

The piecewise-linear test function is constructed as follows:
\begin{align*}
    &V(x):=\max\limits_{\substack{K\subseteq\mathcal K:\\\kappa\in K\Rightarrow p_{\kappa}\in\mathcal K}}\bigg\{\frac{1+(|K|-1)\delta}{|K|}\sum_{k\in K}x_k\bigg\},
\end{align*}
where $\delta$ can be any small value such that $0<\delta<1$.

The following lemmas are useful in proving Theorem \ref{thm_decentralized2} where we consider the regime $X\subseteq\mathcal Q$ containing $x$.
\begin{lmm}
\label{lmm_holding}
A bottleneck can not be in the holding status.
\end{lmm}
\noindent\emph{Proof.} Otherwise, the bottleneck must have at least one downstream subserver. By \eqref{invariant}, $x_{k_s}\geq x_b$. Since $b$ is a bottleneck, we have
\begin{align}\label{eq_bottleneck}
    \frac{1+(|K^X|-2)\delta}{|K^X|-1}\sum_{k\neq b}x_k\leq\frac{1+(|K^X|-1)\delta}{|K^X|}\sum_{k\in K^X}x_k,
\end{align}
which implies
\begin{align}\label{eq_x_k}
    (1-\delta)\sum_{k\in K^X}x_{k}\le|K^X|[1+(|K^X|-2)\delta]x_b.
\end{align}
Since $x_b\le x_{k_s}$, we have
\begin{align*}
    (1-\delta)\sum_{k\in K^X}x_{k}<|K^X|(1+|K^X|\delta)x_{k_s},
\end{align*}
which is equivalent to 
\begin{align*}
    \frac{1+(|K^X|-1)\delta}{|K^X|}\sum_{k\in K^X}x_k<\frac{1+|K^X|\delta}{|K^X|+1}\Big(\sum_{k\in K^X}x_k+x_{k_s}\Big),
\end{align*}
contradicting with the fact that subserver $b$ is dominant and subserver $k_s$ is non-dominant.
\hfill$\square$

\begin{cor}
\label{cor_bottleneck}
Based on \eqref{eq_x_k}, we have $x_b>0$, i.e., any bottleneck $b$ must be non-empty. 

This corollary and Lemma \ref{lmm_holding} ensure that all bottlenecks can discharge customers and contribute negative terms to the drift.
\end{cor}

\begin{lmm}
\label{lmm_non_dominant}
Let $k_r^1$ be the first subserver on route $r$, then either the route with the smallest $x_{k_r^1}$ is non-dominant or every route is dominant.
\end{lmm}

\noindent\emph{Proof.} If there is only one route, then that route must be dominant. Now assume there are at least two routes and route $\hat{r}$ has the smallest $x_{k_r^1}$, i.e. $\forall\ r\in\mathcal R$, $x_{k_{\hat{r}}^1}\leq x_{k_{r}^1}$. Suppose $k_{\hat{r}}^1\in K^X$ and $\exists\ r\in\mathcal R$ s.t. $k_{r}^1\notin K^X$. Note that by \eqref{invariant}, $x_b\leq x_{k_{\hat{r}}^1}\leq x_{k_{r}^1}$, then from \eqref{eq_bottleneck} we have
\begin{align*}
    \frac{1+(|K^X|-1)\delta}{|K^X|}\sum_{k\in K^X}x_k<\frac{1+|K^X|\delta}{|K^X|+1}\Big(\sum_{k\in K^X}x_k+x_{k_{r}^1}\Big),
\end{align*}
contradicting with our supposition.
Therefore, either $\hat{r}$ is non-dominant or every route in $\mathcal R$ is dominant.
\hfill$\square$

\begin{lmm}
\label{lmm_cut}
If $x\in\mathcal{X}$ makes every route $r\in\mathcal R$ dominant, then we have
\begin{align*}
    \sum_{k\in\mathcal B^X}\mu_{k}(x)
    =\sum_{n:n=n_k,k\in\mathcal B^X}\bar\mu_{n}
\end{align*}
\end{lmm}

\noindent\emph{Proof.}
Once there is an inactive bottleneck $k$ and a non-dominant but active duplicate subserver $k'$, the imaginary switch mechanism will move the job being served in $k'$ to $k$ and move one job in $k$ to $k'$. This is allowed since both $k$ and $k'$ contain at least one job due to the fact that a job being served in $k'$ and the bottleneck $k$ must be non-empty.
\hfill$\square$

Similar to the proof of Theorem~\ref{thm_centralized}, we first analyze the internal transmissions and then the external arrivals.

\subsubsection{Internal transmissions}
In the proof of Theorem~\ref{thm_centralized}, we have already discussed the case where internal transmissions between subservers are on the same route. However, unlike the JSR policy, the JSQ-AS policy allows internal transmissions between subservers on different routes. Hence, we also need to consider the internal transmission from subserver $k$ to subserver $j$ where $r_k\neq r_j$. 

The definition of dominance ensures that if $k$ is non-dominant, so is $k_s$. According to the routing policy, $x_{k_s}\geq x_j$. Let $\ell$ be the first non-dominant subserver on route $r_k$ and $b$ be the bottleneck on route $r_j$. If $j$ is non-dominant, then by \eqref{invariant}, we have $x_\ell\geq x_{k_s}\geq x_j\geq x_b$. Now from \eqref{eq_bottleneck} we can obtain
\begin{align*}
    \frac{1+(|K^X|-1)\delta}{|K^X|}\sum_{k\in K^X}x_k<\frac{1+|K^X|\delta}{|K^X|+1}\Big(\sum_{k\in K^X}x_k+x_{\ell}\Big),
\end{align*}
contradicting with the definition of dominant subservers. 

Thus, it cannot be the case that $k$ is non-dominant and $j$ is dominant, which implies that any internal transmission does not positively contribute to the mean drift.

\subsubsection{External arrivals}
According to Lemma \ref{lmm_non_dominant}, if a non-dominant route exists, then the routing policy guarantees that an arriving job must be routed to the first subserver on a non-dominant route $r$; this leads to non-positive contribution to the mean drift. Otherwise, every route is dominant. Then for any $x\in\mathcal Q$ ($x\neq0$), the drift satisfies
\begin{align*}
    D^X(x)
    &\overset{\text{Corollary } \ref{cor_bottleneck}}{\le}\frac{1+(|K^X|-1)\delta}{|K^X|}\Big(\lambda-\sum_{b\in\mathcal B^X}\mu_{b}(x)\zeta_b(x)\Big)\\
    &\overset{\eqref{invariant}}{=}\frac{1+(|K^X|-1)\delta}{|K^X|}\Big(\lambda-\sum_{b\in\mathcal B^X}\mu_{b}(x)\Big)\\
    &\overset{\text{Lemma }\ref{lmm_cut}}{=}\frac{1+(|K^X|-1)\delta}{|K^X|}\Big(\lambda-\sum_{n:n=n_k,k\in\mathcal B^X}\bar\mu_{n}\Big)\\
    &\overset{\text{Lemma }\ref{lmm_stabilizability}}{<}0,
\end{align*}
which completes the proof.
\hfill$\square$

The JSQ-AS policy cannot be directly applied to multi-class network, because the imaginary switch mechanism may move a job to the subserver of a different class with a different destination. Although the imaginary service rate control in the JSR policy can be used for multiple classes, it needs global information such as the information of dominance and bottlenecks for the preemption, so it is not suitable for the decentralized setting. The design of a decentralized MDI control policy for multi-class network can be a future work. 
\section{Concluding remarks}
\label{sec_conclude}

We study the stability of open queuing networks under a class of model data-independent control policies. In addition, we derive an easy-to-use stability criterion based on route expansion of the network and explicit piecewise-linear test functions.
With the stability criterion, we generalize the classical join-the-shortest-queue policy to ensure stability and attain maximum throughput under centralized/decentralized settings.
Our analysis and design can also be applied to specific network control problems with stability issues. 

\bibliographystyle{IEEEtran}  
\bibliography{Bibliography}

\begin{thebibliography}{10}
\providecommand{\url}[1]{#1}
\csname url@samestyle\endcsname
\providecommand{\newblock}{\relax}
\providecommand{\bibinfo}[2]{#2}
\providecommand{\BIBentrySTDinterwordspacing}{\spaceskip=0pt\relax}
\providecommand{\BIBentryALTinterwordstretchfactor}{4}
\providecommand{\BIBentryALTinterwordspacing}{\spaceskip=\fontdimen2\font plus
\BIBentryALTinterwordstretchfactor\fontdimen3\font minus
  \fontdimen4\font\relax}
\providecommand{\BIBforeignlanguage}[2]{{%
\expandafter\ifx\csname l@#1\endcsname\relax
\typeout{** WARNING: IEEEtran.bst: No hyphenation pattern has been}%
\typeout{** loaded for the language `#1'. Using the pattern for}%
\typeout{** the default language instead.}%
\else
\language=\csname l@#1\endcsname
\fi
#2}}
\providecommand{\BIBdecl}{\relax}
\BIBdecl

\bibitem{kumar1995stability}
P.~Kumar and S.~P. Meyn, ``Stability of queueing networks and scheduling
  policies,'' \emph{IEEE Transactions on Automatic Control}, vol.~40, no.~2,
  pp. 251--260, 1995.

\bibitem{meyn2001sequencing}
S.~P. Meyn, ``Sequencing and routing in multiclass queueing networks part i:
  Feedback regulation,'' \emph{SIAM Journal on Control and Optimization},
  vol.~40, no.~3, pp. 741--776, 2001.

\bibitem{smith2010dynamic}
S.~L. Smith, M.~Pavone, F.~Bullo, and E.~Frazzoli, ``Dynamic vehicle routing
  with priority classes of stochastic demands,'' \emph{SIAM Journal on Control
  and Optimization}, vol.~48, no.~5, pp. 3224--3245, 2010.

\bibitem{zhang2018analysis}
R.~Zhang, F.~Rossi, and M.~Pavone, ``Analysis, control, and evaluation of
  mobility-on-demand systems: a queueing-theoretical approach,'' \emph{IEEE
  Transactions on Control of Network Systems}, vol.~6, no.~1, pp. 115--126,
  2018.

\bibitem{bertsimas1994optimization}
D.~Bertsimas, I.~C. Paschalidis, and J.~N. Tsitsiklis, ``Optimization of
  multiclass queueing networks: Polyhedral and nonlinear characterizations of
  achievable performance,'' \emph{The Annals of Applied Probability}, pp.
  43--75, 1994.

\bibitem{down1997piecewise}
D.~Down and S.~P. Meyn, ``Piecewise linear test functions for stability and
  instability of queueing networks,'' \emph{Queueing Systems}, vol.~27, no.
  3-4, pp. 205--226, 1997.

\bibitem{gallager2013stochastic}
R.~G. Gallager, \emph{Stochastic processes: theory for applications}.\hskip 1em
  plus 0.5em minus 0.4em\relax Cambridge University Press, 2013.

\bibitem{dai1995positive}
J.~G. Dai, ``On positive {Harris} recurrence of multiclass queueing networks:
  {A} unified approach via fluid limit models,'' \emph{The Annals of Applied
  Probability}, pp. 49--77, 1995.

\bibitem{foss1998stability}
S.~Foss and N.~Chernova, ``On the stability of a partially accessible
  multi-station queue with state-dependent routing,'' \emph{Queueing Systems},
  vol.~29, no.~1, pp. 55--73, 1998.

\bibitem{tang2020analysis}
Y.~Tang and L.~Jin, ``Analysis and control of dynamic flow networks subject to
  stochastic cyber-physical disruptions,'' \emph{arXiv preprint
  arXiv:2004.00159}, 2020.

\bibitem{sarikaya2011dynamic}
Y.~Sarikaya, T.~Alpcan, and O.~Ercetin, ``Dynamic pricing and queue stability
  in wireless random access games,'' \emph{IEEE Journal of Selected Topics in
  Signal Processing}, vol.~6, no.~2, pp. 140--150, 2011.

\bibitem{dube2009bertrand}
P.~Dube and R.~Jain, ``Bertrand games between multi-class queues,'' in
  \emph{Proceedings of the 48h IEEE Conference on Decision and Control (CDC)
  held jointly with 2009 28th Chinese Control Conference}.\hskip 1em plus 0.5em
  minus 0.4em\relax IEEE, 2009, pp. 8588--8593.

\bibitem{savla2011dynamical}
K.~Savla and E.~Frazzoli, ``A dynamical queue approach to intelligent task
  management for human operators,'' \emph{Proceedings of the IEEE}, vol. 100,
  no.~3, pp. 672--686, 2011.

\bibitem{foschini1978basic}
G.~Foschini and J.~Salz, ``A basic dynamic routing problem and diffusion,''
  \emph{IEEE Transactions on Communications}, vol.~26, no.~3, pp. 320--327,
  1978.

\bibitem{ephremides1980simple}
A.~Ephremides, P.~Varaiya, and J.~Walrand, ``A simple dynamic routing
  problem,'' \emph{IEEE transactions on Automatic Control}, vol.~25, no.~4, pp.
  690--693, 1980.

\bibitem{vvedenskaya1996queueing}
N.~D. Vvedenskaya, R.~L. Dobrushin, and F.~I. Karpelevich, ``Queueing system
  with selection of the shortest of two queues: An asymptotic approach,''
  \emph{Problemy Peredachi Informatsii}, vol.~32, no.~1, pp. 20--34, 1996.

\bibitem{eschenfeldt2018join}
P.~Eschenfeldt and D.~Gamarnik, ``Join the shortest queue with many servers.
  the heavy-traffic asymptotics,'' \emph{Mathematics of Operations Research},
  vol.~43, no.~3, pp. 867--886, 2018.

\bibitem{gupta2007analysis}
V.~Gupta, M.~H. Balter, K.~Sigman, and W.~Whitt, ``Analysis of
  join-the-shortest-queue routing for web server farms,'' \emph{Performance
  Evaluation}, vol.~64, no. 9-12, pp. 1062--1081, 2007.

\bibitem{mukhopadhyay2015analysis}
A.~Mukhopadhyay and R.~R. Mazumdar, ``Analysis of randomized
  join-the-shortest-queue (jsq) schemes in large heterogeneous
  processor-sharing systems,'' \emph{IEEE Transactions on Control of Network
  Systems}, vol.~3, no.~2, pp. 116--126, 2015.

\bibitem{mehdian2017join}
S.~Mehdian, Z.~Zhou, and N.~Bambos, ``Join-the-shortest-queue scheduling with
  delay,'' in \emph{2017 American Control Conference (ACC)}.\hskip 1em plus
  0.5em minus 0.4em\relax IEEE, 2017, pp. 1747--1752.

\bibitem{tang2020security}
Y.~Tang, Y.~Wen, and L.~Jin, ``Security risk analysis of the shorter-queue
  routing policy for two symmetric servers,'' in \emph{2020 American Control
  Conference (ACC)}.\hskip 1em plus 0.5em minus 0.4em\relax IEEE, 2020, pp.
  5090--5095.

\bibitem{bramson2011stability}
M.~Bramson \emph{et~al.}, ``Stability of join the shortest queue networks,''
  \emph{The Annals of Applied Probability}, vol.~21, no.~4, pp. 1568--1625,
  2011.

\bibitem{dai2007stability}
J.~Dai, J.~J. Hasenbein, and B.~Kim, ``Stability of join-the-shortest-queue
  networks,'' \emph{Queueing Systems}, vol.~57, no.~4, pp. 129--145, 2007.

\bibitem{foley2001join}
R.~D. Foley and D.~R. McDonald, ``Join the shortest queue: stability and exact
  asymptotics,'' \emph{The Annals of Applied Probability}, vol.~11, no.~3, pp.
  569--607, 2001.

\bibitem{ling2010global}
X.~Ling, M.-B. Hu, R.~Jiang, and Q.-S. Wu, ``Global dynamic routing for
  scale-free networks,'' \emph{Physical Review E}, vol.~81, no.~1, p. 016113,
  2010.

\bibitem{ren2011traffic}
F.~Ren, T.~He, S.~K. Das, and C.~Lin, ``Traffic-aware dynamic routing to
  alleviate congestion in wireless sensor networks,'' \emph{IEEE Transactions
  on Parallel and Distributed Systems}, vol.~22, no.~9, pp. 1585--1599, 2011.

\bibitem{papadimitriou1994complexity}
C.~H. Papadimitriou and J.~N. Tsitsiklis, ``The complexity of optimal queueing
  network control,'' in \emph{Proceedings of IEEE 9th Annual Conference on
  Structure in Complexity Theory}.\hskip 1em plus 0.5em minus 0.4em\relax IEEE,
  1994, pp. 318--322.

\bibitem{towsley1980queuing}
D.~Towsley, ``Queuing network models with state-dependent routing,''
  \emph{Journal of the ACM (JACM)}, vol.~27, no.~2, pp. 323--337, 1980.

\bibitem{kelly1993dynamic}
F.~Kelly and C.~Laws, ``Dynamic routing in open queueing networks: Brownian
  models, cut constraints and resource pooling,'' \emph{Queueing systems},
  vol.~13, no. 1-3, pp. 47--86, 1993.

\bibitem{sarachik1982decentralized}
P.~Sarachik and U.~Ozguner, ``On decentralized dynamic routing for congested
  traffic networks,'' \emph{IEEE Transactions on Automatic Control}, vol.~27,
  no.~6, pp. 1233--1238, 1982.

\bibitem{gregoire2014capacity}
J.~Gregoire, X.~Qian, E.~Frazzoli, A.~De~La~Fortelle, and T.~Wongpiromsarn,
  ``Capacity-aware backpressure traffic signal control,'' \emph{IEEE
  Transactions on Control of Network Systems}, vol.~2, no.~2, pp. 164--173,
  2014.

\bibitem{varaiya2013max}
P.~Varaiya, ``Max pressure control of a network of signalized intersections,''
  \emph{Transportation Research Part C: Emerging Technologies}, vol.~36, pp.
  177--195, 2013.

\bibitem{dai1995stability}
J.~G. Dai and S.~P. Meyn, ``Stability and convergence of moments for multiclass
  queueing networks via fluid limit models,'' \emph{IEEE Transactions on
  Automatic Control}, vol.~40, no.~11, pp. 1889--1904, 1995.

\end{thebibliography}

\end{document}